\documentclass{svproc}

\usepackage{amssymb,amsmath}

\usepackage{tikz}

\usepackage[]{algorithm2e}

\usepackage[hidelinks]{hyperref}
\usepackage{cleveref}
\crefformat{equation}{\textup{#2(#1)#3}}
\crefrangeformat{equation}{\textup{#3(#1)#4--#5(#2)#6}}
\crefmultiformat{equation}{\textup{#2(#1)#3}}{ and \textup{#2(#1)#3}}
{, \textup{#2(#1)#3}}{, and \textup{#2(#1)#3}}
\crefrangemultiformat{equation}{\textup{#3(#1)#4--#5(#2)#6}}%
{ and \textup{#3(#1)#4--#5(#2)#6}}{, \textup{#3(#1)#4--#5(#2)#6}}{, and \textup{#3(#1)#4--#5(#2)#6}}

\Crefformat{equation}{#2Equation~\textup{(#1)}#3}
\Crefrangeformat{equation}{Equations~\textup{#3(#1)#4--#5(#2)#6}}
\Crefmultiformat{equation}{Equations~\textup{#2(#1)#3}}{ and \textup{#2(#1)#3}}
{, \textup{#2(#1)#3}}{, and \textup{#2(#1)#3}}
\Crefrangemultiformat{equation}{Equations~\textup{#3(#1)#4--#5(#2)#6}}%
{ and \textup{#3(#1)#4--#5(#2)#6}}{, \textup{#3(#1)#4--#5(#2)#6}}{, and \textup{#3(#1)#4--#5(#2)#6}}

\crefname{algocf}{alg.}{algs.}
\Crefname{algocf}{Algorithm}{Algorithms}

\crefname{figure}{Fig.}{Figs.}

\usepackage{subcaption}

\usepackage{url}

\usepackage{color}

\newcommand{\sDomain}{\Omega}
\newcommand{\sbDomain}{{\partial\sDomain}}
\newcommand{\sbDomainIn}{{\partial\sDomain_{\text{in}}}}
\newcommand{\sbDomainOut}{{\partial\sDomain_{\text{out}}}}
\newcommand{\pDomain}{\Gamma}
\newcommand{\pTrainSet}{{\pDomain_{\text{train}}}}

\newcommand{\para}{\xi}
\newcommand{\paraIn}{\para^{\text{in}}}
\newcommand{\paraOut}{\para^{\text{out}}}
\newcommand{\paraOutInd}[1]{\paraOut_{#1}}
\newcommand{\paraRef}{{\xi_{\text{ref}}}}

\newcommand{\soO}{\Theta}
\newcommand{\outcome}{\theta}
\newcommand{\sA}{\mathcal{F}}
\newcommand{\pM}{\mathbb{P}}

\newcommand{\hSpace}{X}
\newcommand{\feSpace}{\hSpace_\feStep}
\newcommand{\feSpaceK}{\feSpace^\tNumb}
\newcommand{\feSpaceKn}{\feSpace^{\tNumb+1}}

\newcommand{\dSpace}{\hSpace '}

\newcommand{\rbSpace}{\hSpace_\rbDim}

\newcommand{\rbSpaceKn}{\rbSpace^{\tNumb+1}}
\newcommand{\rbSpacen}{\hSpace_{\rbDim+1}}
\newcommand{\drbSpace}{\widetilde{\hSpace}_{\drbDim}}

\newcommand{\drbSpaceKn}{\drbSpace^{\tNumb+1}}
\newcommand{\rbSpaceSol}{\hSpace^\sol_\rbDim}
\newcommand{\rbSpaceSolw}{\hSpace_{\rbDim}^{\sol,\pdf}}

\newcommand{\rbSpaceOutP}{\hSpace^\qoi_{\rbDim}}
\newcommand{\rbSpaceOutD}{\widetilde{\hSpace}^\qoi_{\drbDim}}

\newcommand{\rbSpaceOutwP}{\hSpace^{\qoi,\pdf}_{\rbDim}}
\newcommand{\rbSpaceOutwD}{\widetilde{\hSpace}^{\qoi,\pdf}_{\drbDim}}
\newcommand{\podSpace}{\hSpace_{\text{POD},\rbDim}}

\newcommand{\oFunc}{l}
\newcommand{\qoi}{s}
\newcommand{\feQoi}{\qoi_\feStep}
\newcommand{\rbQoi}{\qoi_\rbDim}

\newcommand{\res}{r}

\newcommand{\rbRes}{\res_\rbDim}
\newcommand{\drbRes}{\tilde{\res}_\drbDim}
\newcommand{\fcRes}{\tilde{\res}_{\text{fc}}}
\newcommand{\blf}{a}

\newcommand{\rhsFunc}{b}

\newcommand{\pdf}{\rho}
\newcommand{\feBf}{\phi}				
\newcommand{\rbBf}{\varphi}				
\newcommand{\drbBf}{\zeta}				

\newcommand{\podSol}{\sol_{\text{POD},\rbDim}}					
\newcommand{\podBf}{w}					
			
\newcommand{\podSV}{\lambda}				
		
\newcommand{\klInd}{l}				
\newcommand{\klRf}{\kappa}				
\newcommand{\klEigF}{c}					
\newcommand{\klEigV}{\mu}				
\newcommand{\klExp}{\alpha_0}				
\newcommand{\klCovK}{C}					
\newcommand{\klNumb}{{N_{\text{KL}}}}			

\newcommand{\mcExp}{\EE_{\text{MC}}}			
\newcommand{\mcNumb}{{N_{\text{MC}}}}			

\newcommand{\UnifD}[2]{\mathcal{U}\left(#1,#2\right)}
\newcommand{\BetaD}[4]{\mathcal{B}\left(#1,#2,#3,#4\right)}

\DeclareMathOperator{\spn}{span}
\DeclareMathOperator*{\argmax}{arg\,max}
\newcommand\Sqrt[1]{#1^{1/2}}

\newcommand\SSqrt[1]{\left(#1\right)^{1/2}}
\newcommand{\abs}[1]{\lvert#1\rvert}
\newcommand{\tDer}{\partial_\tVar}

\newcommand{\xVar}{x}
\newcommand{\feStep}{h}
\newcommand{\tVar}{t}
\newcommand{\tEnd}{T}
\newcommand{\tInd}{k}
\newcommand{\tNumb}{K}

\newcommand{\tStep}{\Delta\tVar}

\newcommand{\sol}{u}
\newcommand{\feSol}{\sol_{\feStep}}
\newcommand{\rbSol}{\sol_\rbDim}
\newcommand{\wrbSol}{\sol_{\rbDim,\pdf}}
\newcommand{\dSol}{\psi}
\newcommand{\dfeSol}{\dSol_\feStep}
\newcommand{\drbSol}{\dSol_\drbDim}
\newcommand{\feSolC}[1]{\sol_{\feStep,#1}}	
\newcommand{\rbSolC}[1]{\sol_{\rbDim,#1}}	
\newcommand{\dfeSolC}[1]{\dSol_{\feStep,#1}}	
\newcommand{\drbSolC}[1]{\dSol_{\drbDim,#1}}	

\newcommand{\coerConst}{\alpha}
\newcommand{\contConst}{\gamma}
\newcommand{\lbCoerConst}{\overline{\coerConst}}
\newcommand{\ubContConst}{\overline{\contConst}}
\newcommand{\eps}{\epsilon}
\newcommand{\epsTol}{\eps_{\text{tol}}}

\newcommand{\pDim}{p}
\newcommand{\feDim}{{\mathcal{N}}}
\newcommand{\rbDim}{N}
\newcommand{\drbDim}{{\tilde{\rbDim}}}

\newcommand{\RR}{\mathbb R}
\newcommand{\NN}{\mathbb N}
\newcommand{\EE}{\mathbb E}

\newcommand{\vertiii}[1]{{\left\vert\kern-0.25ex\left\vert\kern-0.25ex\left\vert #1 
\right\vert\kern-0.25ex\right\vert\kern-0.25ex\right\vert}}
\newcommand\norm[1]{\left\lVert#1\right\rVert}
\newcommand\hNorm[1]{\norm{#1}_{\hSpace}}
\newcommand\dNorm[1]{\norm{#1}_{\hSpace'}}

\newcommand\eNormPara[1]{\vertiii{#1}_\para^{\text{pr}}}
\newcommand\eNormParaSingle[1]{\norm{#1}_\para}
\newcommand\eNormParaSingleRef[1]{\norm{#1}_\paraRef}
\newcommand\deNormPara[1]{\vertiii{#1}_\para^{\text{du}}}

\newcommand\LtwoNorm[1]{\norm{#1}_{L^2(\sDomain)}}

\newcommand{\Ltwo}[2]{(#1,#2)_{L^2(\sDomain)}}

\newcommand{\iProd}[2]{(#1,#2)}

\newcommand{\iparaProd}[2]{\iProd{#1}{#2}_\para}
\newcommand{\iparaRefProd}[2]{\iProd{#1}{#2}_\paraRef}

\newcommand{\est}{\Delta}
\newcommand{\estDim}{\Delta_\rbDim}
\newcommand{\estSol}{\est_\rbDim^{\sol}}
\newcommand{\estSolw}{\est_\rbDim^{\sol,\pdf}}
\newcommand{\estdSol}{\est_\drbDim^{\dSol}}
\newcommand{\estOut}{\est_{\rbDim,\drbDim}^{\qoi}}
\newcommand{\estOutw}{\est_{\rbDim,\drbDim}^{\qoi,\pdf}}
\newcommand{\estFc}{\est_\drbDim^{\dSol,\text{fc}}}

\begin{document}
\mainmatter              
\title{A weighted reduced basis method for parabolic PDEs
with random data
}
\titlerunning{A weighted reduced basis method}  
%
\author{Christopher Spannring\inst{1,2}, Sebastian Ullmann\inst{1,2}  \and
Jens Lang\inst{1,2}}
\authorrunning{Christopher Spannring et al.} 
%
%
\institute{Graduate School of Computational Engineering, Technische Universit\"at Darmstadt,
Dolivostr. 15, 62493, Darmstadt, Germany \\
spannring@gsc.tu-darmstadt.de \\
WWW home page: http://www.graduate-school-ce.de/index.php?id=688
\and
Department of Mathematics, Technische Universit\"at Darmstadt,
Dolivostr. 15, 62493, Darmstadt, Germany}

\maketitle              

\begin{abstract}
This work considers a weighted POD-greedy method to estimate statistical outputs parabolic PDE problems with parametrized random data.
The key idea of weighted reduced basis methods is to weight the parameter-dependent error estimate according to a probability measure in the set-up of the reduced space.
The error of stochastic finite element solutions is usually measured in a root mean square sense regarding their dependence on the stochastic input parameters.
An orthogonal projection of a snapshot set onto a corresponding POD basis defines an optimum reduced approximation in terms of a Monte Carlo discretization of the root mean square error.
The errors of a weighted POD-greedy Galerkin solution are compared against an orthogonal projection of the underlying snapshots onto a POD basis for a numerical example involving thermal conduction.
In particular, it is assessed whether a weighted POD-greedy solutions is able to come significantly closer to the optimum than a non-weighted equivalent.
Additionally, the performance of a weighted POD-greedy Galerkin solution is considered with respect to the mean absolute error of an adjoint-corrected functional of the reduced solution.
\end{abstract}
\section{Introduction}
Because the computational complexity of numerical simulations is growing, model order reduction becomes an essential task.
In the last decades the reduced basis method (RBM) was extensively developed, 
e.g. recent overviews can be found in \cite{HesthavenRozzaStamm2016,QuarteroniManzoniNegri2016,Haadonk2017}.
Furthermore, RBM was applied to stochastically influenced parametrized partial differential equations (PPDEs). 
A review can be found in \cite{ChenQuarteroniRozza2017}.
In \cite{ChenQuarteroniRozza2013} the idea of a weighted RBM for elliptic PDEs was introduced. It allows to build up
more efficient reduced spaces regarding an approximation of statistical quantities. In the following the focus is on the approximation of the 
expectation. The question arises, how much faster does the expected error converge for 
a weighted approach and how close is the reduced solution to an optimal reduced solution.
Using a proper orthogonal decomposition (POD) \cite{KunischVolkwein2002}, an optimal reduced space concerning the mean square error can be constructed.
Furthermore, by a standard primal-dual approach \cite{PrudhommeRovasVeroy2001} also the expected error of a linear output functional is considered.
This work shows different ways to construct reduced order models (ROMs), where the focus is the approximation of the expected value.

The work is organized as follows. In \cref{sec:plpm} the model problem is formulated, and the
notation for the high dimensional discretization is introduced. In \cref{sec:rbm} the reduced model and its assumptions are described.
Non-weighted and weighted error estimators for the reduced model are stated in \cref{sec:errest}. Further, the reduced space construction
for the non-weighted and the weighted approach is described in \cref{sec:rsc}. 
\Cref{sec:romc} describes a ROM, obtained by a Galerkin projection onto a POD.
Numerical results for an instationary heat equation are presented in \cref{sec:numex}. 

\section{Parametrized Linear Parabolic Model}
\label{sec:plpm}
In the following the weak formulation of a parametrized, linear parabolic PDE is considered, where 
a $\pDim$-dimensional parameter vector is random. Let $(\soO,\sA,\pM)$ be a
complete probability space where the sample space $\soO$ contains all possible outcomes $\outcome\in\soO$.
The sigma algebra $\sA$ is given by a subset of all possible subsets of $\soO$, i.e. $\sA\subseteq 2^\soO$,
and the probability measure $\pM\colon\sA\to[0,1]$ maps an event to its probability. 
Let $\para\colon\soO\to\pDomain$
denote a random parameter vector whose image lies in a given parameter domain $\pDomain\subset\RR^\pDim$,
which is determined by the support of the random variables. We assume that $\para$ has a joint probability
density function (pdf) $\pdf\colon\pDomain\to\RR^+$ such that $\int_\soO\mathrm{d}\pM(\outcome)=\int_\pDomain\pdf(\para)\mathrm{d}\para=1$.

The interest lies not only in the solution of the parabolic PDE problem itself but rather in some parameter-dependent output 
$\qoi\colon\pDomain\to\RR$. 
It is computed by an linear (output) functional $\oFunc\colon\hSpace\to\RR$ that maps the solution, lying in a Sobolev space $\hSpace$,
to a scalar output. The continuous parametrized problem reads as follows:
For given $\para\in\pDomain$, compute
\begin{equation*}
  \qoi(\para)=\oFunc(\sol(\tEnd;\para)),
\end{equation*}
where for any $\tVar\in[0,\tEnd]$ the solution $\sol(\tVar;\para)\in\hSpace$
fulfills
\begin{alignat*}{3}
   \langle\partial_\tVar\sol(\tVar;\para),v\rangle+\blf(\sol(\tVar;\para),v)&=\rhsFunc(v;\para), &\quad\forall v\in\hSpace, \\
  \Ltwo{\sol(0;\para)}{v} &=0, &\quad\forall v\in\hSpace.
\end{alignat*}
Here, $\langle\cdot,\cdot\rangle$ denotes a duality pairing between $\dSpace$ and $\hSpace$. The time derivative of the solution needs to lie in 
the dual space, 
i.e. $\partial_\tVar\sol(\tVar;\para)\in \dSpace$ for all parameters.
In order to guarantee existence and uniqueness of the problem,
the bilinear form is uniformly coercive and uniformly bounded 
and the functionals are uniformly bounded, i.e.
\begin{alignat}{3}
  &\lbCoerConst\leq\coerConst(\para)=\inf_{0\neq v\in\hSpace}\frac{\blf(v,v;\para)}{\hNorm{v}^2},\quad &&
  \sup_{0\neq\sol,v\in\hSpace}\frac{\blf(\sol,v;\para)}{\hNorm{\sol}\hNorm{v}}=\contConst(\para)\leq\ubContConst,\text{ and} \label{eq:coercontConst}\\
  &\sup_{0\neq v\in\hSpace}\frac{|\rhsFunc(v;\para)|}{\hNorm{v}}=\contConst_\rhsFunc(\para)\leq\ubContConst_\rhsFunc,\quad && 
  \sup_{0\neq v\in\hSpace}\frac{|\oFunc(v)|}{\hNorm{v}}=\contConst_\oFunc\leq\ubContConst_\oFunc ,
\end{alignat}
with $0<\lbCoerConst$ and $\ubContConst,\ubContConst_\rhsFunc,\ubContConst_\oFunc<\infty,$ 
for all $\para\in\pDomain$ . Furthermore, the bilinear form is assumed to be symmetric, i.e. $\blf(u,v;\cdot)=\blf(v,u;\cdot),\,\forall u,v\in\hSpace$.
It defines a parameter-dependent and parameter-independent energy norm, such that
\begin{alignat}{2}
    \eNormParaSingle{v}&=\sqrt{\iparaProd{v}{v}}=\sqrt{\blf(v,v;\para)},\quad&&\forall v\in\hSpace,\,\forall\para\in\pDomain, \label{eq:eNormParaSingle} \\
    \eNormParaSingleRef{v}&=\sqrt{\iparaRefProd{v}{v}}=\sqrt{\blf(v,v;\paraRef)},\quad&&\forall v\in\hSpace \label{eq:eNormParaSingleRef}.
\end{alignat}
The parameter-independent norm is determined by a fixed reference parameter $\paraRef\in\pDomain$.

\subsection{Fully Discretized Problem}
\label{sec:fdp}
The problem is discretized with an implicit Euler method in time and with a linear finite element method in space.
The time interval $[0,\tEnd]$ is split into $\tNumb\in\NN$ equidistant time intervals with time step size 
$\tStep:=\frac{\tEnd}{\tNumb}$. The solution is approximated at time points 
$\{\tVar^\tInd=\tInd\tStep\colon\tInd=0,\dots,\tNumb\}$, such that $\sol(\xVar,\tVar^\tInd;\para)\approx\sol^\tInd(\xVar;\para)$.
For the space discretization, $\feStep$ denotes the spatial step size and the high dimensional space
$\feSpace:=\spn\{\feBf_1,\dots,\feBf_\feDim\}\subset\hSpace$, with $\dim(\feSpace)=\feDim$, 
contains piecewise linear basis functions. Further, $\feSpaceKn:=\feSpace\times\dots\times\feSpace$ 
denotes the $(\tNumb+1)$th power of the discretized space.
The fully discretized primal problem reads as follows: For given $\para\in\pDomain$, compute 
\begin{equation}
  \label{eq:fdProbOut}
  \feQoi(\para)=\oFunc(\feSol^\tNumb(\para)),
\end{equation}
where the detailed solutions $\{\feSol^\tInd(\para)\}_{\tInd=0}^\tNumb\in\feSpaceKn$ fulfill
\begin{alignat}{3}
  \label{eq:fdProba}
    \Ltwo{\feSol^\tInd(\para)}{v}&+\tStep\,\blf(\feSol^\tInd(\para),v;\para) &&\nonumber\\
    &=\Ltwo{\feSol^{\tInd-1}(\para)}{v}+\tStep\,\rhsFunc(v;\para),\quad&&\forall v\in\feSpace,\,\tInd=1,\dots,\tNumb,\\
    \Ltwo{\feSol^\tInd(\para)}{v}&=0, \quad && \forall v\in\feSpace,\,\tInd=0. \label{eq:fdProbb}
\end{alignat}
A finite element method entails an $\feDim$-dimensional system of linear algebraic equations and hence the discretized problem is computationally 
expensive to solve. The solution coefficients, coming out of the algebraic equations, uniquely represent the detailed solution, such that
\begin{equation*}
\feSol^\tInd(\para)=\sum_{i=1}^\feDim\feSolC{i}^\tInd(\para)\feBf_i,\quad \tInd=0,\dots,\tNumb.
\end{equation*}
\Cref{eq:fdProbOut,eq:fdProba,eq:fdProbb} are referred to as the detailed 
model and $\feSol^\tInd(\cdot)$ as the detailed solution. It can be seen as a reference solution. In the following
the error is measured between the solution of the detailed model and the reduced model.

In order to achieve higher accuracy for the output computation, a dual (or adjoint) problem is used \cite[Chapter 2.1]{GreplPatera2005}.
This is a standard approach in the context of error analysis for functionals \cite{OdenPrudhomme2001}.
The dual problem of \cref{eq:fdProbOut,eq:fdProba,eq:fdProbb} reads as follows: For given $\para\in\pDomain$,
find the dual solutions $\{\dfeSol^\tInd(\para)\}_{\tInd=0}^{\tNumb}\in\feSpaceKn$, s.t.
\begin{alignat}{3}
 \Ltwo{v}{\dfeSol^{\tInd}(\para)}&+\tStep\>\blf(v,\dfeSol^{\tInd}(\para);\para) &&\nonumber\\
			&=\Ltwo{v}{\dfeSol^{\tInd+1}(\para)}, \quad &&\forall v\in\feSpace,\,\tInd=0,\dots,\tNumb-1, \label{eq:fdProbc}\\
			\Ltwo{v}{\dfeSol^{\tInd}(\para)}&=\oFunc(v), \quad&&\forall v\in\feSpace,\,\tInd=\tNumb. \label{eq:fdfinCond}
\end{alignat}
Due to the parameter-independent functional $\oFunc(\cdot)$, the solution of the final condition in \cref{eq:fdfinCond} is parameter-independent, i.e.
$\dfeSol^{\tNumb}(\para)=\dfeSol^{\tNumb}$.
Note, the dual problem evolves backward in time, hence the solutions are computed for decreasing time index $\tInd$.
As for the primal problem, the dual solution has a unique representation
\begin{equation*}
\dfeSol^\tInd(\para)=\sum_{i=1}^\feDim\dfeSolC{i}^\tInd(\para)\feBf_i,\quad \tInd=0,\dots,\tNumb. 
\end{equation*}

\section{Reduced Basis Method}
\label{sec:rbm}
In order to decrease the computation time, a ROM is sought. This can be achieved
using the RBM \cite{HesthavenRozzaStamm2016,QuarteroniManzoniNegri2016}. It
is a Galerkin projection onto a reduced basis space $\rbSpace:=\spn\{\rbBf_1,\dots,\rbBf_\rbDim\}\subset\feSpace$. 
The basis functions are orthogonal w.r.t. \cref{eq:eNormParaSingleRef} and the dimension
is denoted by $\dim(\rbSpace)=\rbDim$. 
Additionally, a reduced space for the dual problem $\drbSpace:=\spn\{\drbBf_1,\dots,\drbBf_\drbDim\}\subset\feSpace$ is introduced.
The basis functions are orthogonal w.r.t. \cref{eq:eNormParaSingleRef} and the dimension
is denoted by $\dim(\drbSpace)=\drbDim$.  
The $(\tNumb+1)$th power of the reduced spaces are defined by
$\rbSpaceKn:=\rbSpace\times\dots\times\rbSpace$ and $\drbSpaceKn:=\drbSpace\times\dots\times\drbSpace$ respectively.
Then, the reduced primal problem reads as follows: For given $\para\in\pDomain$,  
find the reduced solutions $\{\rbSol^\tInd(\para)\}_{\tInd=0}^\tNumb\in\rbSpaceKn$, s.t.
\begin{alignat}{3}
    \Ltwo{\rbSol^\tInd(\para)}{v}&+\tStep\,\blf(\rbSol^\tInd(\para),v;\para) &&\nonumber\\
    &=\Ltwo{\rbSol^{\tInd-1}(\para)}{v}+\tStep\,\rhsFunc(v;\para),~&&\forall v\in\rbSpace,\,\tInd=1,\dots,\tNumb,  \label{eq:rProba}\\
    \Ltwo{\rbSol^\tInd(\para)}{v}&=0, \,&& \forall v\in\rbSpace,~\tInd=0. \label{eq:rProbb}
\end{alignat}
The reduced dual problem reads as follows: For given $\para\in\pDomain$,  
find the reduced solutions
$\{\drbSol^\tInd(\para)\}_{\tInd=0}^{\tNumb}\in\drbSpaceKn$, s.t.
\begin{alignat}{3}
 \Ltwo{v}{\drbSol^{\tInd}(\para)}&+\tStep\>\blf(v,\drbSol^{\tInd}(\para);\para) && \nonumber\\
    &=\Ltwo{v}{\drbSol^{\tInd+1}(\para)},~ &&\forall v\in\drbSpace,\,\tInd=0,\dots,\tNumb-1, \label{eq:rProbc}\\
			\Ltwo{v}{\drbSol^{\tInd}}&=\oFunc(v),~&&\forall v\in\drbSpace,\,\tInd=\tNumb \label{eq:rfinCond}.
\end{alignat}
Note, the reduced spaces $\rbSpace$ and $\drbSpace$ are spanned by a different basis and in general 
they can have different
dimensions, i.e. $\rbDim\neq\drbDim$.
The reduced primal problem and the reduced dual problem yield a system of linear equations with dimensions 
$\rbDim$ and $\drbDim$ respectively. 
For each time step and each parameter the solutions of the algebraic equations yield the primal solution coefficients $\{\rbSolC{n}^\tInd(\para)\}_{n=1}^\rbDim$ and the dual solution coefficients $\{\drbSolC{n}^\tInd(\para)\}_{n=1}^\drbDim$. They determine a unique representation of the reduced solutions, for $\tInd=0,\dots,\tNumb$ and $\para\in\pDomain$, such that
\begin{alignat*}{3}
\rbSol^\tInd(\para)=\sum_{n=1}^\rbDim\rbSolC{n}^\tInd(\para)\rbBf_n,\quad
\drbSol^\tInd(\para)=\sum_{n=1}^\drbDim\drbSolC{n}^\tInd(\para)\drbBf_n.  
\end{alignat*}

The residual for the reduced primal problem is defined by, for $\tInd=1,\dots,\tNumb$,
\begin{alignat*}{2}
    \rbRes^\tInd(v;\para)&=\rhsFunc(v;\para)-\frac{1}{\tStep}\Ltwo{\rbSol^\tInd(\para)-\rbSol^{\tInd-1}(\para)}{v}-\blf(\rbSol^\tInd(\para),v;\para),~&&\forall v\in\feSpace.\\
    \intertext{The residual for the reduced dual problem is defined by, for $\tInd=0,\dots,\tNumb-1$,}
    \drbRes^{\tInd}(v;\para)&=-\frac{1}{\tStep}\Ltwo{v}{\drbSol^{\tInd}(\para)-\drbSol^{\tInd+1}(\para)}-\blf(v,\drbSol^{\tInd}(\para);\para),~&&\forall v\in\feSpace.
\end{alignat*}
Note, the primal residual and the dual residual are orthogonal onto their reduced spaces, i.e. $\rbSpace\subset\ker(\rbRes^\tInd(\cdot;\para))$ 
and $\drbSpace\subset\ker(\drbRes^\tInd(\cdot;\para))$. In addition, the residual of the final condition \cref{eq:rfinCond} is given by,
\begin{equation}
\label{eq:rfinCondRes}
 \fcRes(v)=\oFunc(v)-\Ltwo{v}{\drbSol^{\tNumb}}, \quad \forall v\in\feSpace.
\end{equation}
The reduced output is determined by
\begin{equation}
\label{eq:rProbOut}
  \rbQoi(\para)=\oFunc(\rbSol^\tNumb(\para))+\tStep\,\sum_{\tInd=1}^\tNumb\rbRes^\tInd(\drbSol^{\tInd-1}(\para);\para).
\end{equation}
The second term in \cref{eq:rProbOut} is a ``correction term'', that uses the reduced dual solutions
$\{\drbSol^\tInd\}_{k=0}^{\tNumb-1}$ in order to
achieve higher accuracy for the output computation. Typically the ``correction'' doubles the order of
accuracy for the output approximation, cf. \cite{PierceGiles2000}.

For computational efficiency, the computation is split into an offline phase and an online phase.
The former is expensive to compute and depends
  on the large dimension $\feDim$, and it is related to the reduced model construction. 
Once the reduced model exists, solutions are obtained very fast in the online phase by calculations depending only on the 
  reduced dimension $\rbDim$. 
  For such a splitting, the bilinear form $\blf(\cdot,\cdot;\para)$ and 
  the functional $\rhsFunc(\cdot;\para)$ need to be affine with respect to $\para$ (also known as parameter separable), i.e.
     \begin{alignat}{2} 
   \label{eq:blfps}
      \blf(v,w;\para)&=\sum_{q=1}^{Q_\blf}\theta_q^\blf(\para)\blf_q(v,w),\quad&&\forall v,w\in\hSpace,\,\forall\para\in\pDomain,\\
   \label{eq:rhsps}
      \rhsFunc(v;\para)&=\sum_{q=1}^{Q_\rhsFunc}\theta_q^\rhsFunc(\para)\rhsFunc_q(v),\quad&&\forall v\in\hSpace,\,\forall\para\in\pDomain.   
    \end{alignat}
If this assumption does not hold, an empirical interpolation method (EIM) \cite{BarraultMadayNguyenPatera2004} can be used instead.
    
\section{Error Estimation}
\label{sec:errest}
The objective in this section is the a posteriori error estimation in order to assess the accuracy of the reduced models.
Meaning, the error between the detailed
output and the reduced output is measured in a given norm.
Rigorous error bounds for the solution error and the output error will be stated. 
By the assumptions for an offline-online decomposition in \cref{eq:blfps,eq:rhsps}, the error bounds are computationally inexpensive,
compared to the exact error computation.
The error bounds are computed by the dual norm
of the residuals.
The following error estimators can be extended for a non-symmetric bilinear form $\blf(\cdot,\cdot;\para)$.

\subsection{Non-Weighted Error Estimators}
\label{sec:nwerrest}
The following statements are taken from \cite{GreplPatera2005}:
For the final condition \cref{eq:rfinCond}, the error can be estimated by,
\begin{equation}
  \label{eq:finCondEst}
 \LtwoNorm{\dfeSol^{\tNumb}-\drbSol^{\tNumb}}\leq\estFc:=\sup_{v\in\feSpace}\frac{|\fcRes(v)|}{\LtwoNorm{v}},
\end{equation}
where the final condition residual \cref{eq:rfinCondRes} is maximized. 
Note, if the reduced space for the dual problem contains $\dfeSol^{\tNumb}\in\drbSpace$, the error estimator of the final condition is zero, 
i.e. $\estFc=0$.
The solution errors of the primal problem and the dual problem are estimated in parameter-dependent energy norms, 
\begin{align*}
      \eNormPara{v}&:=\SSqrt{\LtwoNorm{v^\tNumb}^2+\tStep\,\sum_{\tInd=1}^{\tNumb}\eNormParaSingle{v^{\tInd}}^2},\quad&&\forall v\in\feSpaceK, \\
      \deNormPara{v}&:=\SSqrt{\LtwoNorm{v^0}^2+\tStep\,\sum_{\tInd=0}^{\tNumb-1}\eNormParaSingle{v^{\tInd}}^2},\quad&&\forall v\in\feSpaceK,
\end{align*}
where $\eNormParaSingle{\cdot}$ is defined in \cref{eq:eNormParaSingle}.
Rigorous error estimators for the primal solution, dual solution and output for all $\para\in\pDomain$, are defined by,
\begin{align}
 \eNormPara{\feSol(\para)-\rbSol(\para)}&\leq\estSol(\para):=\Sqrt{\left(\frac{\tStep}{\coerConst(\para)}\sum_{\tInd=1}^\tNumb\dNorm{\rbRes^\tInd(\cdot;\para)}^2\right)}, \label{eq:errest}\\
 \deNormPara{\dfeSol(\para)-\drbSol(\para)}&\leq\estdSol(\para):=\Sqrt{\left(\frac{\tStep}{\coerConst(\para)}\sum_{\tInd=0}^{\tNumb-1}\dNorm{\drbRes^\tInd(\cdot;\para)}^2+(\estFc)^2\right)}, \label{eq:derrest}\\
 \left|\feQoi(\para)-\rbQoi(\para)\right|&\leq\estOut(\para):=\estSol(\para)\estdSol(\para). \label{eq:outest}
\end{align}
For the error estimators the coercivity constant $\alpha(\para)$, 
defined in \cref{eq:coercontConst}, comes in. 
It can be approximated by a successive constraint method (SCM) \cite{HuynhRozzaSenPatera2007} for instance.
The dual norms of the residuals in \cref{eq:errest,eq:derrest} are computed by means of the Riesz representation
theorem (see, e.g., \cite[Chapter 2.4]{BrennerScott2008}).

\subsection{Weighted Error Estimators}
\label{sec:werrest}
The pdf $\pdf(\cdot)$ appears in the computation of statistical quantities. Using the results from the previous section, 
the expected solution error can be estimated by,
\begin{equation*}
 \EE[\eNormPara{\feSol-\rbSol}]=\int_\pDomain\eNormPara{\feSol(\para)-\rbSol(\para)}\pdf(\para)\mathrm{d}\para\leq\int_\pDomain\estSol(\para)\pdf(\para)\mathrm{d}\para,
\end{equation*}
just as the expected output error
\begin{equation*}
 \EE[|\feQoi-\rbQoi|]=\int_\pDomain|\feQoi(\para)-\rbQoi(\para)|\pdf(\para)\mathrm{d}\para\leq\int_\pDomain\estOut(\para)\pdf(\para)\mathrm{d}\para.
\end{equation*}
The following weighted error estimators, introduced in \cite{ChenQuarteroniRozza2013}, are defined,
\begin{alignat}{2}
 \estSolw(\para)&:=\estSol(\para)\pdf(\para), &&\quad\forall\para\in\pDomain, \label{eq:werrest}\\
 \estOutw(\para)&:=\estOut(\para)\pdf(\para), &&\quad\forall\para\in\pDomain, \label{eq:woutest}
\end{alignat}
where the weight $\pdf(\para)$ gives greater weight to more likely parameter values. 
The weighted estimators \cref{eq:werrest,eq:woutest} are still computationally cheap and will be used as a optimality criterion 
for a weighted reduced space construction.

\section{Reduced Space Construction}
\label{sec:rsc}
In this section, a non-weighted and a weighted reduced space construction are described.  
The reduced space construction is based on the POD-greedy algorithm \cite{HaasdonkOhlberger2008}, stated in
\cref{alg:podg}. 
\begin{algorithm}[]
  \DontPrintSemicolon
   \KwData{$\epsTol$, parameter training set $\pTrainSet\subset\pDomain$, $\para^{(1)}$}
  \KwResult{reduced space $\rbSpace$}
  $\rbDim=1,\,\hSpace_1=\spn\{\feSol^\tNumb(\para^{(1)})\}$\;
  \While{
  $\eps_\rbDim:=\max_{\para\in\pTrainSet}\estDim(\para)\,>\epsTol$}{ 
  \BlankLine
    $\para^{(\rbDim+1)}:=\argmax_{\para\in\pTrainSet}\estDim(\para)$ \; \BlankLine  
  $\text{compute }\feSol^\tInd(\para^{(\rbDim+1)}),~ k=0,\dots,\tNumb$, using \cref{eq:rProba,eq:rProbb}  \; \BlankLine  
  $e_P^\tInd(\para^{(\rbDim+1)}):=\feSol^\tInd(\para^{(\rbDim+1)})-P_{\rbSpace}\feSol^\tInd(\para^{(\rbDim+1)}),~ \tInd=0,\dots,\tNumb$ \; \BlankLine  
  $\rbBf_{\rbDim+1}:=POD_1(\{e_P^\tInd(\para^{(\rbDim+1)})\}_{\tInd=0}^\tNumb)$  \; \BlankLine  
 $\rbSpacen:=\rbSpace\oplus\spn\{\rbBf_{\rbDim+1}\}$\; \BlankLine
 $N:=N+1$\; \BlankLine}
  \caption{POD-greedy algorithm}
  \label{alg:podg}
 \end{algorithm}
 
Note, the algorithm is formulated for the primal problem \cref{eq:fdProba,eq:fdProbb,eq:rProba,eq:rProbb}. 
Analogously the algorithm can be utilized for the dual problem \cref{eq:fdProbc,eq:fdfinCond,eq:rProbc,eq:rfinCond}. 
 
The POD-greedy algorithm combines a greedy algorithm for the parameter domain and a POD for the time interval. 
It computes the detailed solution $\feSol^\tNumb(\para^{(1)})$, for a
given parameter $\para^{(1)}$, which spans the initial reduced space $\hSpace_1$.
The dimension of the reduced basis $\rbDim$ grows iteratively.
As a stopping criterion for the iteration, an error estimator of the previous section 
needs to fall below some given error tolerance $\epsTol>0$.
Alternatively, a predefined reduced dimension can be given as a stopping criterion.
The parameter value $\para^{\rbDim+1}$ is determined in each iteration by evaluating an optimality criterion. 
Maximizing the exact errors over a large parameter domain can be computationally infeasible. 
Instead, the 
computationally cheap error estimators are maximized, which is often called weak POD-greedy. 
The maximum is sought over a training parameter set, which is a
finite and uniformly sampled approximation set of the parameter domain.
The solutions for the single time steps $\{\feSol^\tInd(\para^{(\rbDim+1)})\}_{\tInd=1}^\tNumb$, called snapshots, 
are computed.
In order to compress the information of the obtained solution trajectory, 
the first POD mode of the projection error
onto the reduced space is computed and added to the reduced basis. 

\subsection{A Non-weighted Reduced Space Construction}
\label{sec:nwrsc}
A non-weighted reduced space construction is based on \cref{alg:podg}.
It is distinguished if the objective is either the approximation of the solution
$\feSol$ or the approximation of the output $\feQoi$.
The former uses the primal solution error estimator in \cref{eq:errest} for the optimality criterion in the POD-greedy procedure.
The latter uses the output error estimator in \cref{eq:outest} for the optimality criterion in the POD-greedy procedure.
The output error estimator consists of a primal error estimator \cref{eq:errest} and dual error estimator \cref{eq:derrest}. Therefore, maximizing the output
error estimator yields a reduced space regarding the primal problem 
and a reduced space regarding the dual problem.
Note, a non-weighted reduced space construction weights all parameter values $\para\in\pDomain$ equally.

\subsection{A Weighted Reduced Space Construction}
\label{sec:wrsc}
As in the previous section, a weighted reduced space construction is based on \cref{alg:podg}. However, in this section the objective is to build up a reduced space, that gives better error convergence rates 
regarding statistical 
quantities, compared to the non-weighted approach.
Since the input parameters are random, certain parameters are more likely to appear. 
Highly probable parameters obtain more importance
incorporating the pdf. Hence, weighted error estimators of \cref{sec:werrest}
are used for the reduced space construction.
It is distinguished if the objective is either the approximation of the expected solution
$\EE[\feSol]$ or the approximation of the expected output $\EE[\feQoi]$.
The former uses the weighted primal solution error estimator in \cref{eq:werrest} for the optimality criterion in the POD-greedy procedure.
The latter uses the weighted output error estimator in \cref{eq:woutest} for the optimality criterion in the POD-greedy procedure.
The weighted output error estimator consists of a primal error estimator \cref{eq:errest} and dual error estimator \cref{eq:derrest}. Therefore, maximizing the weighted output
error estimator yields a weighted reduced space regarding the primal problem 
and a weighted reduced space regarding the dual problem.
Note, the weighted approach uses the same uniformly sampled parameter set $\pTrainSet$ as the non-weighted approach.
Meaning, the weighting comes in only by the optimality criterion that maximizes the weighted error estimators.

\section{ROM Comparison With A POD Projection}
\label{sec:romc}
In this work, a conceptual comparison between different ROMs is drawn.
In \cref{sec:nwrsc,sec:wrsc}
a non-weighted and a weighted ROM construction by a POD-greedy approach were shown. 
In this section, the idea is to build up another reduced basis that 
yields an optimal reduced solution regarding the expected solution error, namely 
\begin{align}
\label{eq:minPODcont}
\begin{split}
 &\min_{\podBf_1,\dots,\podBf_\rbDim\in\feSpace}
\min_{\substack{\sol^1,\dots,\sol^\tNumb\colon\\ \pDomain\to\spn\{\podBf_1,\dots,\podBf_\rbDim\}}}
\EE\left[\tStep\sum_{\tInd=1}^\tNumb\eNormParaSingleRef{\feSol^\tInd-\sol^\tInd}^2\right],\\
&\text{subject to }\iparaRefProd{\podBf_m}{\podBf_n}=\delta_{mn}.
\end{split}
\end{align}
The norm $\eNormParaSingleRef{\cdot}$ and the inner product $\iparaRefProd{\cdot}{\cdot}$ are defined by \cref{eq:eNormParaSingleRef} and $\delta_{mn}$ is the Kronecker delta. 
An optimal choice for the unknown functions $\{\sol^\tInd\}_{\tInd=1}^\tNumb$ can be achieved 
by choosing their orthogonal projection onto
the
$\rbDim$-dimensional 
POD space $\podSpace:=\spn\{\podBf_1,\dots,\podBf_\rbDim\}$, namely
\begin{equation}
\label{eq:podRepr}
\podSol^\tInd(\para):=\sum_{n=1}^\rbDim\iparaRefProd{\feSol^\tInd(\para)}{\podBf_n}\podBf_n\in\podSpace,\quad \tInd=1,\dots\,\tNumb. 
\end{equation}
The orthogonal projection yields an equivalent formulation of the minimization problem \cref{eq:minPODcont}, such that
\begin{align}
\label{eq:minPOD}
\begin{split}
 &\min_{\podBf_1,\dots,\podBf_\rbDim\in\feSpace}
\EE\left[\tStep\sum_{\tInd=1}^\tNumb\eNormParaSingleRef{\feSol^\tInd-\sum_{n=1}^\rbDim\iparaRefProd{\feSol^\tInd}{\podBf_n}\podBf_n}^2\right],\\
&\text{subject to }\iparaRefProd{\podBf_m}{\podBf_n}=\delta_{mn}.
\end{split}
\end{align}
The expectation in \cref{eq:minPODcont} and \cref{eq:minPOD} can not be determined analytically. Hence, 
the expectation is approximated by a Monte Carlo (MC) method
(see, e.g., \cite[Chapter 2.7.3]{Fishman1996}). It uses snapshots $\{\feSol^\tInd(\para^{(i)})\}_{i=1}^\mcNumb$, with 
random realizations $\{\para^{(i)}\}_{i=1}^\mcNumb$ sampled by its pdf, such that
\begin{equation}
\label{eq:MCest}
 \mcExp[\feSol^\tInd]=\frac{1}{\mcNumb}\sum_{i=1}^{\mcNumb}\feSol^\tInd(\para^{(i)}),\quad\tInd=1,\dots,\tNumb.
\end{equation}
Note, the number of MC samples $\mcNumb$ is assumed to be large enough, such that the MC error can be neglected, 
i.e. $\EE[\cdot]\approx\mcExp[\cdot]$.
By the MC approximation of the expectation the minimization problem in \cref{eq:minPOD} is given by
\begin{equation}
\label{eq:minPODdiscr}
\begin{split}
 &\min_{\podBf_1,\dots,\podBf_\rbDim\in\feSpace}
 \frac{\tStep}{\mcNumb}\sum_{i=1}^{\mcNumb}\sum_{\tInd=1}^{\tNumb}\eNormParaSingleRef{\feSol^\tInd(\para^{(i)})-\sum_{n=1}^{\rbDim}\iparaRefProd{\feSol^\tInd(\para^{(i)})}{\podBf_n}\podBf_n}^2,\\ 
 &\text{subject to }\iparaRefProd{\podBf_m}{\podBf_n}=\delta_{mn}.
\end{split}
 \end{equation}
This minimization problem can be solved by an eigenvalue problem, see, e.g. \cite[Theorem 1.8]{GubischVolkwein2017}.
The resulting eigenfunctions $\{\podBf_n\}_{n=1}^{\rbDim}$ are the first $\rbDim\in\{1,\dots,\mcNumb\tNumb\}$ orthogonal POD modes 
and the eigenvalues $\{\podSV_l\}_{l=1}^{\mcNumb\tNumb}$ determine the error in \cref{eq:minPODdiscr}, such that 
\begin{equation*}
 \frac{\tStep}{\mcNumb}\sum_{i=1}^{\mcNumb}\sum_{\tInd=1}^{\tNumb}\eNormParaSingleRef{\feSol^\tInd(\para^{(i)})-\sum_{n=1}^{\rbDim}\iparaRefProd{\feSol^\tInd(\para^{(i)})}{\podBf_n}\podBf_n}^2
 =\sum_{l=\rbDim+1}^{\mcNumb\tNumb}\podSV_l.
\end{equation*}
The first $\rbDim$ POD modes span the POD space $\podSpace=\spn\{\podBf_1,\dots,\podBf_\rbDim\}$.

However, the optimality pays its price, since the POD needs the detailed solution for each parameter in the parameter domain. Compared to the 
POD-greedy in \cref{sec:rsc}, the evaluation of the error estimator for the parameter values are
computationally cheap and eventually only $\rbDim$ detailed solutions for the reduced space construction
need to be evaluated. 
Note, the comparison only can be done for the expected solution error. For 
the expected output error $\EE[\abs{\feQoi-\rbQoi}]$ an optimal reduced space cannot be found with a POD.

In this section optimality for the mean square error could be stated. Therefore,
the POD method determines an optimal reduced POD space based on a finite set of snapshots and a parameter-independent energy norm. 
Moreover, the optimal representation of the reduced solution
is given by the orthogonal projection onto the reduced POD space. 
The minimal mean square error can be easily computed
by the truncated eigenvalues coming out of the POD method.

\section{Numerical Example}
\label{sec:numex}
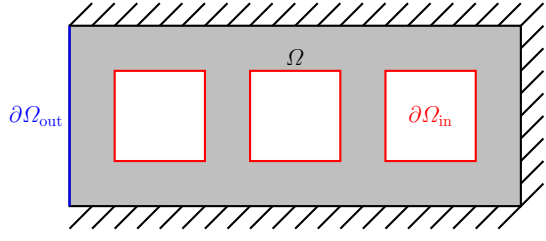
\begin{figure}[t]
\begin{center}
\begin{tikzpicture}[thick,scale=0.6, every node/.style={scale=0.6}]
 \node[rectangle,minimum width=10cm,minimum height=4cm,draw=black,fill=lightgray,label=west:\Large{$\textcolor{blue}{\sbDomainOut}$}] at (0,0) {} ;
 \draw[color=blue] (-5,-2) -- (-5,2) ;
 \node[rectangle,minimum width=2cm,minimum height=2cm,draw=red,fill=white] at (-3,0) {} ;
 \node[rectangle,minimum width=2cm,minimum height=2cm,draw=red,fill=white, label=above:\Large{$\sDomain$}] at (0,0) {} ;
 \node[rectangle,minimum width=2cm,minimum height=2cm,draw=red,fill=white] at (3,0) {\Large{$\textcolor{red}{\sbDomainIn}$}} ;
 \foreach \i in {0,0.5,...,10} 
  \draw (-5+\i,2) -- (-4.5+\i,2.5) ; 
  \foreach \i in {0,0.5,...,9.5}{
  \draw (-4.5+\i,-2) -- (-5+\i,-2.5) ;}
  \foreach \i in {0,0.5,...,3.5}{
  \draw (5,1.5-\i) -- (5.5,2-\i) ; }  
\end{tikzpicture}
\caption{Heat flow in domain $\sDomain$ is considered}
\label{fig:domain}
\end{center}
\end{figure}

An instationary heat transfer is considered \cite{GreplPatera2005} as a numerical example. 
The time dependent heat flow is computed on a rectangular domain where three squares are cut out, see \cref{fig:domain}.
The spatial domain is defined by $\sDomain=\{[0,10]\times[0,4]\}\setminus\{\{(1,3)\times(1,3)\}\cup\{(4,6)\times(1,3)\}\cup\{(7,9)\times(1,3)\}\}$.
As a quantity of interest $\qoi(\cdot)$, the average temperature in the domain at the end time point $\tEnd$ is computed.
The stochastic parameters $\para=(\paraOut,\paraIn)$ enter the problem via the 
boundary condition.
On the left domain boundary $\sbDomainOut$ there is a random heat inflow $\klRf(\cdot;\paraOut)$ modeled by a Karhunen-Lo\`{e}ve (KL) expansion \cite[Chapter 37.5]{Loeve1978}.
Therefore, it needs to hold that $\klRf$ is a second-order random field, meaning that its second moment is finite. As a consequence of the KL,
the corresponding random variables $\paraOut$ are uncorrelated and have zero mean.
The top, right and bottom boundary are insulated. 
The condition at the inner boundary of the domain (squares) $\sbDomainIn$ is parametrized by a beta distributed parameter $\paraIn$, that can be interpreted as a cooling parameter.

\begin{figure}[]
  \begin{center}
  \includegraphics[scale=0.4]{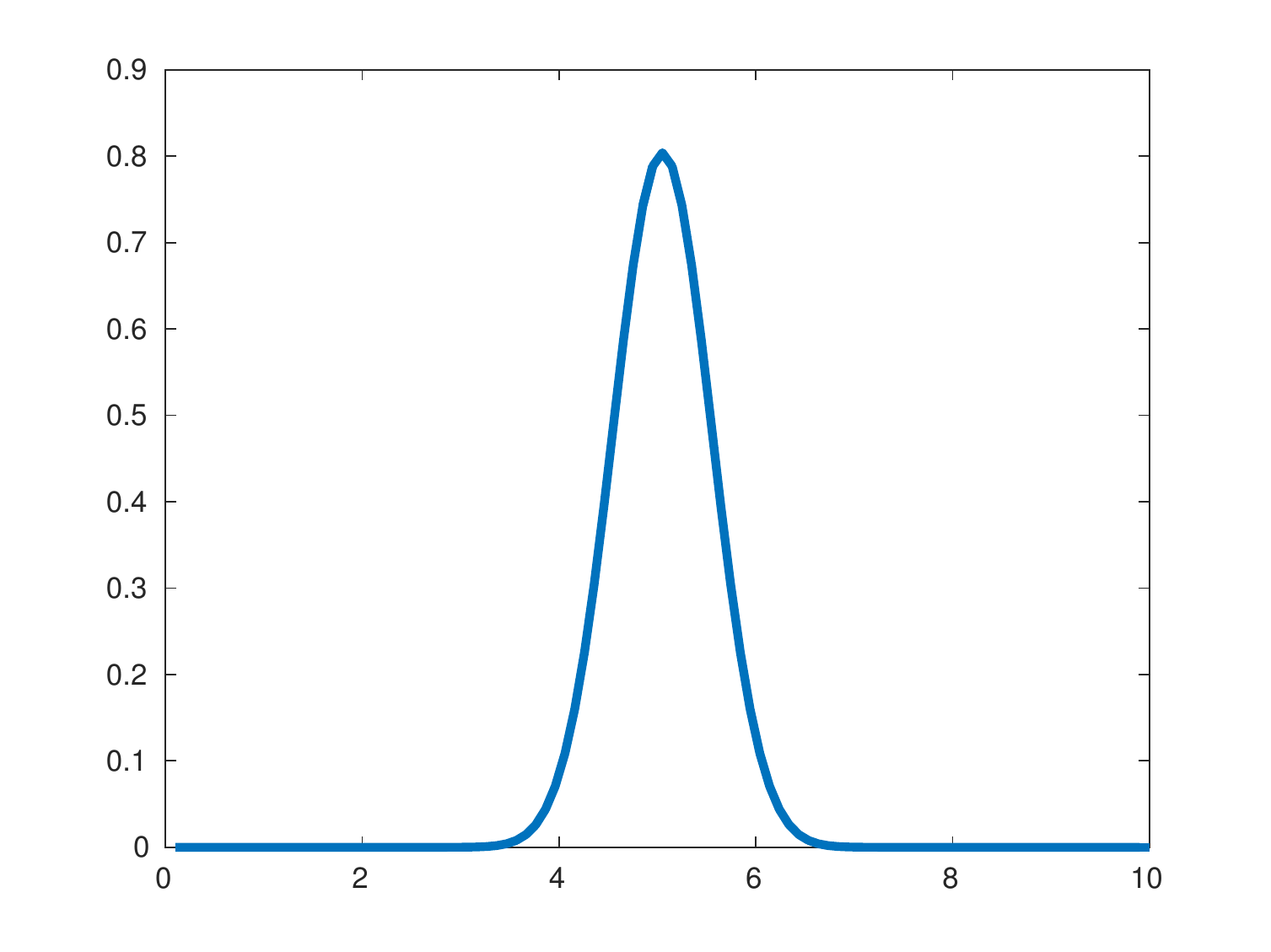} 
  \caption{Pdf of beta distribution with support $[0.1,10]$ and scaling parameter $50$}
  \label{fig:pdf50}
  \end{center}
\end{figure}

The weak formulation of the parabolic problem reads as follows: 
For given realization $\para=\left(\paraOut,\paraIn\right)\in\pDomain$, evaluate 
\begin{equation*}
  \qoi(\para)=\frac{1}{|\sDomain|}\int_\sDomain\sol(\tEnd;\para), 
\end{equation*}
where for any $\tVar\in[0,\tEnd]$ the solution $\sol(t;\para)\in \hSpace:=H^1(\sDomain)$ fulfills
\begin{equation*}
\begin{split}
\int_\sDomain\tDer \sol v + \int_\sDomain\nabla\sol\cdot\nabla v + \int_\sbDomainOut \klRf(\paraOut)\sol v + \paraIn\int_\sbDomainIn\sol v &= \int_\sbDomainOut\klRf(\paraOut) v,\\
\sol(0;\para)&=0,
\end{split}
\end{equation*}
for all $v\in\hSpace$. The boundary integrals 
result from
a parametrized Robin boundary condition,
\begin{equation*}
\frac{\partial\sol}{\partial n}=\begin{cases} 
      \klRf(\paraOut)(1-\sol) & \text{on }\sbDomainOut,\\
			-\paraIn\sol & \text{on }\sbDomainIn, \\
      0 & \text{on }\sbDomain\setminus\{\sbDomainOut\cup\sbDomainIn\}.
   \end{cases}
\end{equation*}
The time interval is determined by $\tEnd=20$, and
the KL field is given by its truncated version $\klRf(x,\paraOut)=\klExp(x)+\sum_{\klInd=1}^{\klNumb}\sqrt{\klEigV_\klInd}\klEigF_\klInd(x)\paraOutInd{\klInd} $ where the first $\klNumb=10$ eigenpairs are chosen.
The eigenvalues $\{\klEigV_l\}_{l=1}^\klNumb$ and the eigenfunctions $\{\klEigF_l\}_{l=1}^\klNumb$ are solutions of an integral equation containing a 
covariance operator \cite[Chapter 7.4]{LordPowellShardlow2014}, corresponding to an exponential kernel 
$\klCovK(x,y)=\exp(-\abs{x-y}/a)$ with correlation length $a=2$.
The expectation of the random field is defined by $\klExp=10$. 
The random parameters $\paraOut=\{\paraOutInd{\klInd}\}_{l=1}^\klNumb$ 
are assumed to be independent of each other and are distributed with 
uniform distribution $\UnifD{-\sqrt{3}}{\sqrt{3}}$.
The random parameter $\paraIn$ is distributed with beta distribution $\BetaD{0.1}{10}{50}{50}$ 
where the first two inputs are the interval bounds of the support,
the third and fourth input are scaling parameters. In \cref{fig:pdf50} its probability density function is drawn.
Further, the parameters $\paraIn$ and $\paraOut$ are independent of each other as well.

\begin{figure}[]
\centering
\begin{subfigure}[]{.5\textwidth}
  \centering
  \includegraphics[width=.85\linewidth]{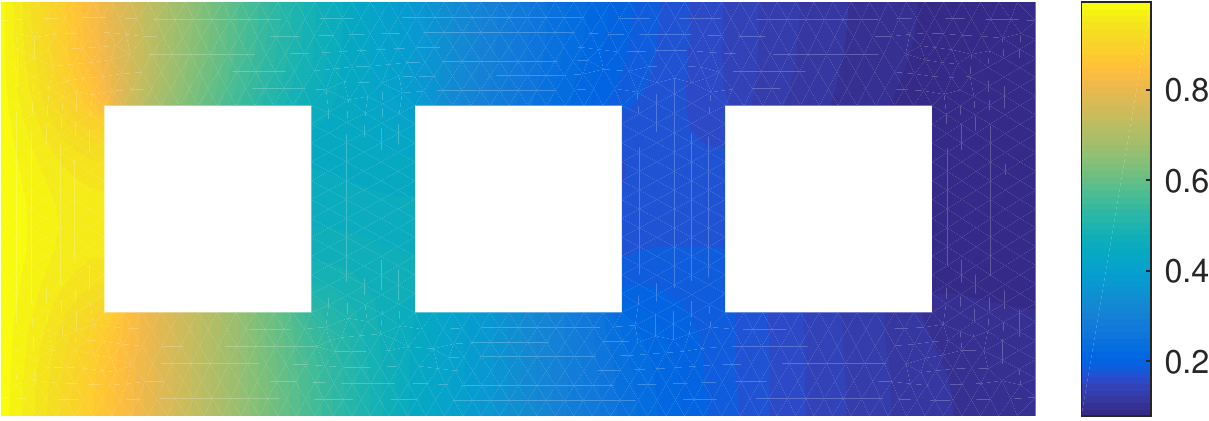}
  \caption{Reference solution}
  \label{fig:refSol}
\end{subfigure}%
\begin{subfigure}[]{.5\textwidth}
  \centering
  \includegraphics[width=.85\linewidth]{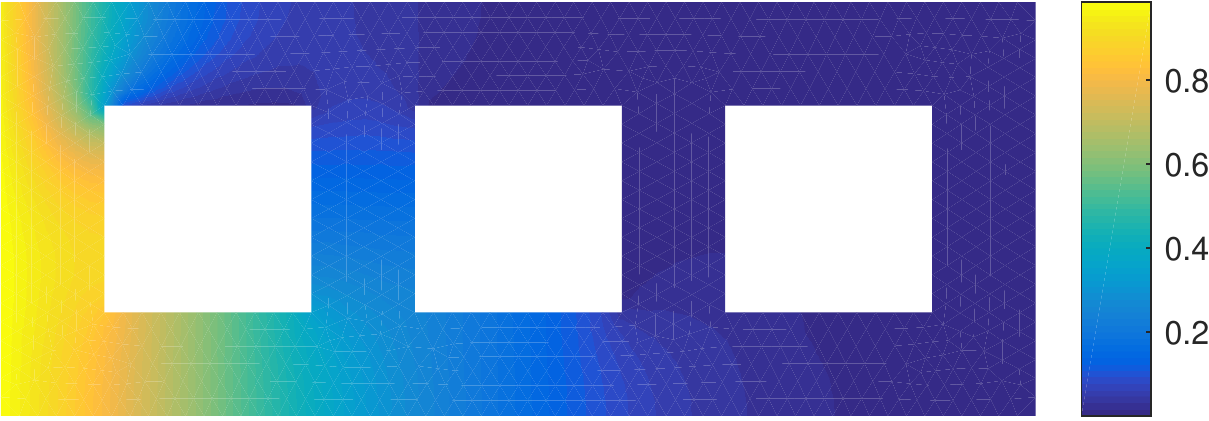}
  \caption{Random solution}
  \label{fig:randSol}
\end{subfigure}
\caption{Finite element solutions at time point $\tEnd$}
\label{fig:refrandsol}
\end{figure}

The bilinear form, right hand side and linear output functional are defined by,
\begin{alignat*}{3}
  \blf(w,v;\para)&:=\int_\sDomain\nabla w\cdot\nabla v + \int_\sbDomainOut \klRf(\paraOut)w v + \paraIn\int_\sbDomainIn w v,\quad &&\forall w,v\in\hSpace,\\
   \rhsFunc(v;\para)&:=\int_\sbDomainOut\klRf(\paraOut) v,\quad\oFunc(v):=\frac{1}{|\sDomain|}\int_\sDomain v,\quad &&\forall v\in\hSpace.\\
   \end{alignat*}
For the spatial discretization, a linear finite element method with $\feDim=1132$ degrees of 
freedom is used. For the time discretization, an implicit Euler method 
is applied with time step size $\tStep=0.2$ and $\tNumb=100$. 
In \cref{fig:refrandsol}, two solutions can be seen, 
the finite element solution for the reference parameter $\paraRef=(0,\dots,0,0.1)\in\pDomain$
in \cref{fig:refSol} and the finite element solution for a randomly sampled parameter in \cref{fig:randSol}.
Recall, the reference parameter defines an inner product, see \cref{eq:eNormParaSingleRef}, such that 
\begin{equation*}
  \iparaRefProd{w}{v}=\int_\sDomain\nabla w\cdot\nabla v + 10 \int_\sbDomainOut wv + 0.1\int_\sbDomainIn w v,\quad \forall w,v\in\hSpace.
\end{equation*}
The 
coercivity constant in \cref{eq:coercontConst} is chosen as a uniformly lower bound, i.e. $\lbCoerConst=1$.
For the reduced space construction, see \cref{sec:rsc},
the parameter domain $\pDomain=\left[-\sqrt{3},\sqrt{3}\right]^\klNumb\times\left[0.1,10\right]$ is approximated by a subset $\pTrainSet\subset\pDomain$, which contains $|\pTrainSet|=500$ 
independent uniformly distributed parameter samples.
All the reduced spaces regarding the primal problem have the same initial basis.
As an initial parameter value
$\para^{(1)}$, the parameter is chosen such
that $\paraIn$ attains its maximum in $\pTrainSet$. The solution at the end time point for that
parameter spans the initial reduced basis, e.g. $\hSpace_1=\spn\{\feSol^\tNumb(\para^{(1)})\}$. 
Further, all the reduced spaces regarding the dual problem have the same initial basis. It 
simply takes the solution of the final condition as initial basis, e.g. $\widetilde{\hSpace}_1=\spn\{\dfeSol^{\tNumb}\}$.
This choice implies that the error estimator for the final condition in \cref{eq:finCondEst} is zero.

In the following different ROMs are compared. Therefore, the root mean square error and the mean absolute output error 
are considered.

\begin{figure}[t]
\centering
\begin{subfigure}{.5\textwidth}
  \centering
  \includegraphics[width=.85\linewidth]{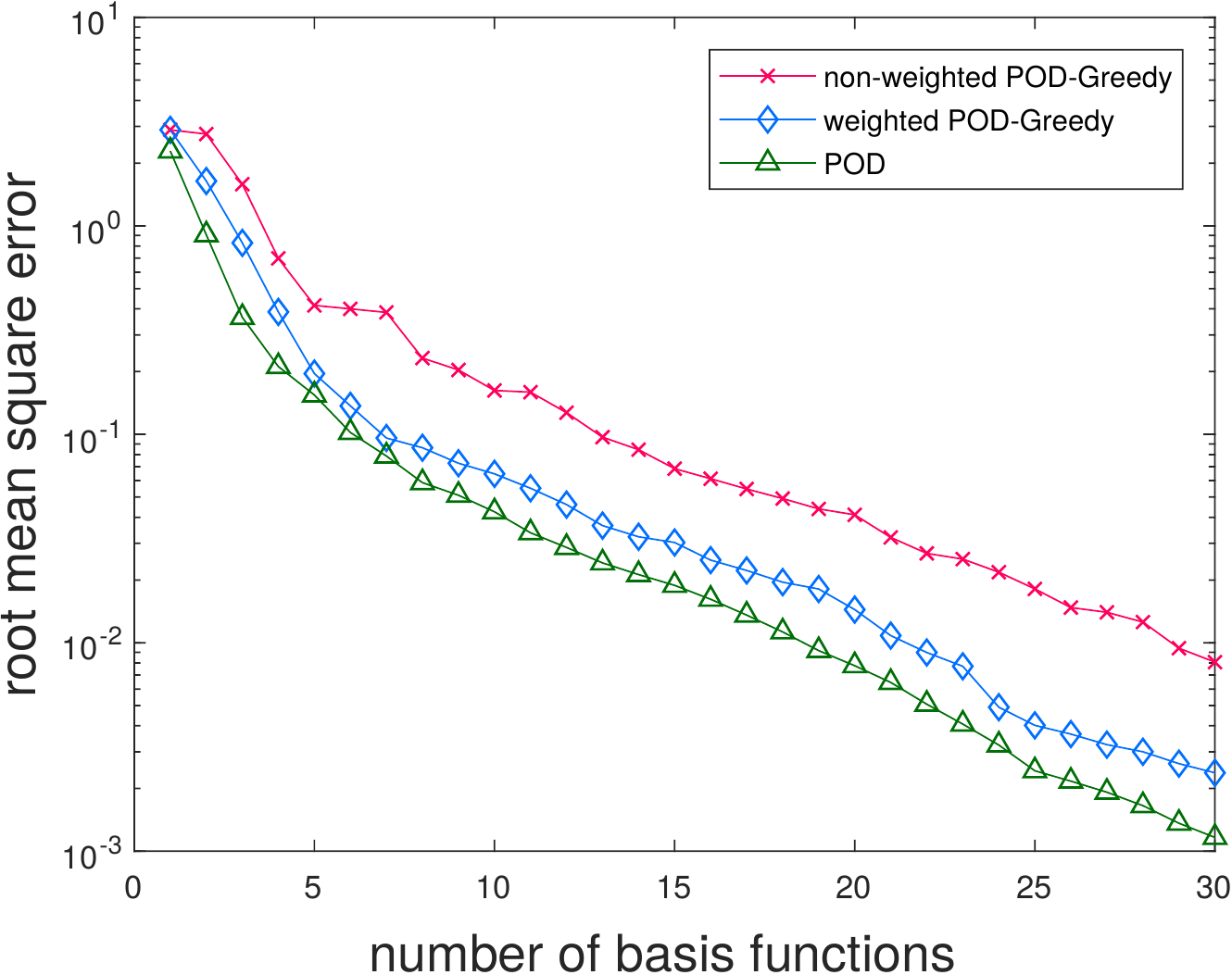}
  \caption{$\Sqrt{\mcExp\left[\tStep\sum_{\tInd=1}^{\tNumb}\eNormParaSingleRef{\feSol^\tInd-\rbSol^\tInd}^2\right]}$}
  \label{fig:expsolerr}
\end{subfigure}%
\begin{subfigure}{.5\textwidth}
  \centering
  \includegraphics[width=.85\linewidth]{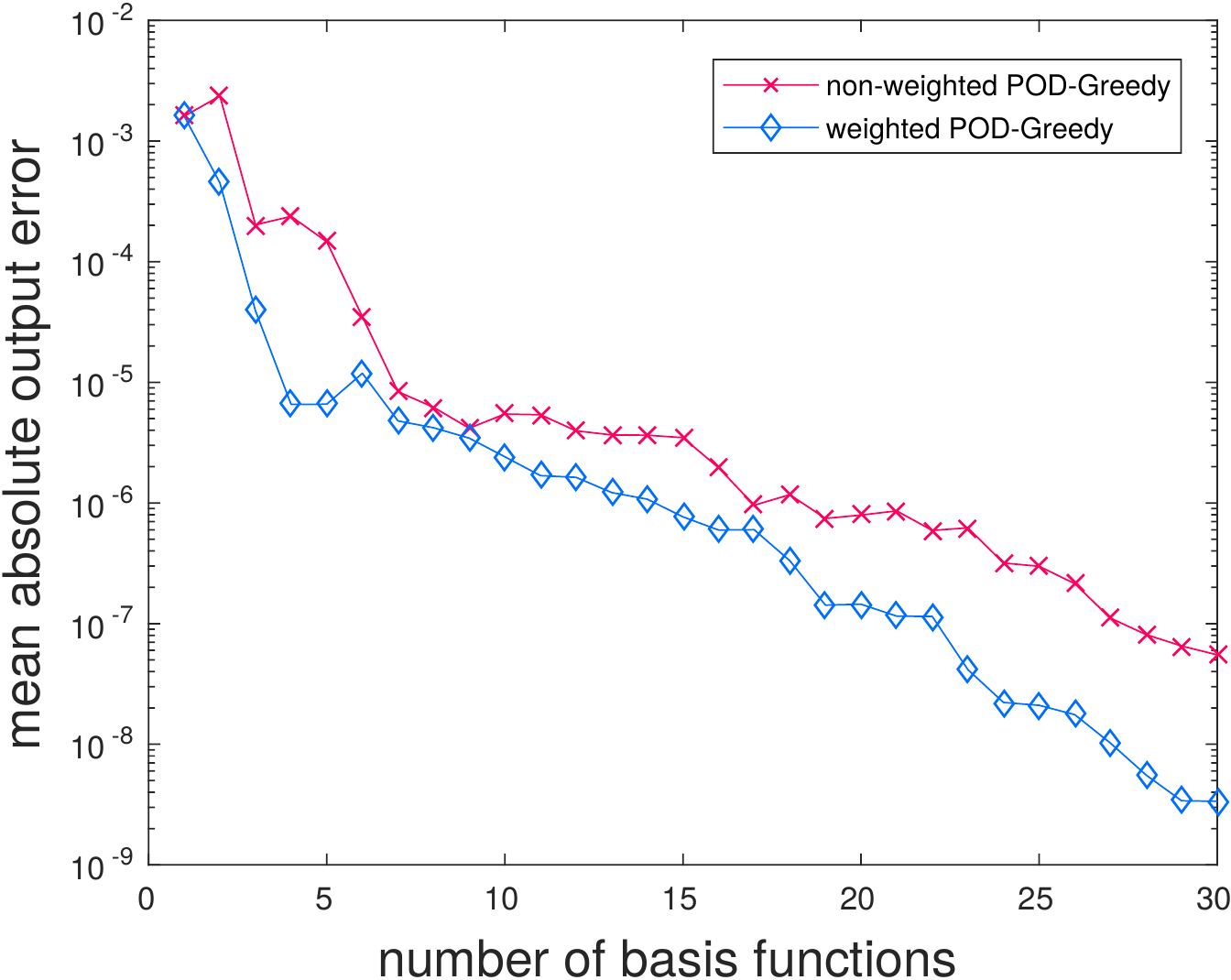}
  \caption{$\mcExp[\abs{\feQoi-\rbQoi}]$}
  \label{fig:expouterr}
\end{subfigure}
\caption{Comparison between non-weighted, weighted POD-greedy and POD}
\label{fig:experrcomp}
\end{figure}

First, the goal is an efficient approximation of the root mean square solution error.
Therefore, three different models are considered. 
The first reduced space $\rbSpaceSol$, is obtained
by a non-weighted POD-greedy algorithm, see \cref{alg:podg}. 
It uses a non-weighted error estimator \cref{eq:errest} as the optimality criterion.
The non-weighted reduced solutions $\{\rbSol^\tInd\in\rbSpaceSol\}_{\tInd=1}^\tNumb$ 
are determined by \cref{eq:rProba,eq:rProbb}.
The second reduced space $\rbSpaceSolw$, is obtained
by a weighted POD-greedy algorithm. 
It uses a weighted error estimator \cref{eq:werrest} as the optimality criterion in \cref{alg:podg}. 
The weighted reduced solutions $\{\wrbSol^\tInd\in\rbSpaceSolw\}_{\tInd=1}^\tNumb$
are determined by \cref{eq:rProba,eq:rProbb}.
The third model is obtained by a POD, see \cref{sec:romc}.
Therefore, the snapshots $\{\feSol^\tInd(\para^{(i)})\}_{i,k=1}^{\mcNumb,\tNumb}$ for each time step and for all parameter values,
sampled by the joint pdf, are computed. The POD basis functions determine
the $POD$ solutions $\{\podSol^\tInd\in\podSpace\}_{\tInd=1}^\tNumb$, see \cref{eq:podRepr}.
Those three solutions determine the expected solution errors in \cref{fig:expsolerr},
where the errors for the first $30$ basis functions are shown.
For the approximation of the expected value, the MC method \cref{eq:MCest} is used. 
Therefore, the same samples as for the POD snapshot space are chosen.

Second, the goal is an efficient approximation of the mean absolute output error, 
where a non-weighted and a weighted approach are compared. 
The former utilizes a non-weighted POD-greedy algorithm, that uses a non-weighted output error estimator 
\cref{eq:outest} as the optimality criterion. 
The output error estimator consists of an error estimator for the primal solution \cref{eq:errest}
and an error estimator for the dual solution \cref{eq:derrest}. 
As explained at the end of \cref{sec:nwrsc}, the non-weighted POD-greedy for the output approximation 
yields two reduced spaces $\rbSpaceOutP$ and $\rbSpaceOutD$.
The non-weighted reduced output is determined by \cref{eq:rProbOut}. The calculation 
of the reduced output is based on the reduced solutions 
$\rbSol^\tNumb\in\rbSpaceOutP$, computed by \cref{eq:rProba,eq:rProbb},
and
$\{\drbSol^\tInd\in\rbSpaceOutD\}_{\tInd=0}^{\tNumb-1}$, computed by \cref{eq:rProbc,eq:rfinCond}. 
The weighted approach utilizes a weighted POD-greedy algorithm, that uses a weighted output error estimator 
\cref{eq:woutest} as the optimality criterion. 
The weighted output error estimator consists of an error estimator for the primal solution \cref{eq:errest}
and an error estimator for the dual solution \cref{eq:derrest}. 
As explained at the end of \cref{sec:wrsc}, the weighted POD-greedy for the output approximation 
yields two reduced spaces $\rbSpaceOutwP$ and $\rbSpaceOutwD$.
The weighted reduced output is determined by \cref{eq:rProbOut}. The calculation 
of the reduced output is based on the reduced solutions 
$\rbSol^\tNumb\in\rbSpaceOutwP$, computed by \cref{eq:rProba,eq:rProbb},
and
$\{\drbSol^\tInd\in\rbSpaceOutwD\}_{\tInd=0}^{\tNumb-1}$, computed by \cref{eq:rProbc,eq:rfinCond}. 
Those two approaches determine the mean absolute output errors in \cref{fig:expouterr},
where the errors for the first $30$ basis functions are shown.
\Cref{fig:experrcomp} shows that in each iteration of the POD-greedy the weighted approach gives better error results compared to the non-weighted approach. Further, the optimal error convergence, obtained from the POD, is observed in \cref{fig:expsolerr}.

\section{Conclusions}
In this work, different model order reduction techniques for a specific parabolic model problem with data uncertainties were studied.
Namely, a RBM, a weighted RBM and a POD were used in order 
to reduce the dimension of a high-fidelity model obtained from a linear finite element model. Apart from the solution, a linear functional maps the solution to a quantity of interest.
The work considered the root mean square 
error and the mean absolute output error. 
A numerical example for an instationary heat transfer with random input data was studied.
It has been shown, that a weighted RBM yields better error results of the root mean square error and 
the mean absolute output error compared to a non-weighted RBM.
The POD yields an optimal reduced space of the mean square error for an energy norm.
The errors for the numerical example show, 
that the weighted RBM is closer to the POD than the non-weighted RBM.\\

\textbf{Acknowledgements. } This work is supported by the Excellence Initiative of the German federal and state governments and the Graduate School of Computational Engineering at the Technische Universit\"at Darmstadt.

%
%


\begin{thebibliography}{6}
%

\bibitem{HesthavenRozzaStamm2016}
Hesthaven, J.S., Rozza, G., Stamm, B.: Certified Reduced Basis Methods for Parametrized Partial Differential Equations.
Springer International Publishing (2016). \url{doi:10.1007/978-3-319-22470-1}

\bibitem{QuarteroniManzoniNegri2016}
Quarteroni, A., Manzoni, A., Negri, F.: Reduced Basis Methods for Partial Differential Equations An Introduction.
Springer International Publishing (2016). \url{doi:10.1007/978-3-319-15431-2}

\bibitem{PierceGiles2000}
Pierce, N.A., Giles, M.B.: Adjoint Recovery of Superconvergent Functionals from PDE Approximations.
SIAM Review (2000). \url{doi:10.1137/S0036144598349423}

\bibitem{OdenPrudhomme2001}
Oden, J.T., Prudhomme, S.: Goal-oriented error estimation and adaptivity for the finite element method.
Computers \& Mathematics with Applications. \url{doi:10.1016/S0898-1221(00)00317-5}

\bibitem{GreplPatera2005}
Grepl, M.A., Patera, A.T.: A posteriori error bounds for reduced-basis approximations of parametrized parabolic partial differential equations.
ESAIM: M2AN (2005). \url{doi:10.1051/m2an:2005006} 

\bibitem{ChenQuarteroniRozza2013}
Chen, P., Quarteroni, A., Rozza, G.: A Weighted Reduced Basis Method for Elliptic Partial Differential Equations with Random Input Data.
SIAM Journal on Numerical Analysis (2013). \url{doi:10.1137/130905253} 

\bibitem{HaasdonkOhlberger2008}
Haasdonk, B., Ohlberger, M.: Reduced basis method for finite volume approximations of parametrized linear evolution equations.
ESAIM: M2AN (2008). \url{doi:10.1051/m2an:2008001} 

\bibitem{KunischVolkwein2002}
Kunisch, K., Volkwein, S.: Galerkin Proper Orthogonal Decomposition Methods for a General Equation in Fluid Dynamics.
SIAM Journal on Numerical Analysis (2002). \url{doi:10.1137/S0036142900382612} 

\bibitem{Loeve1978}
Lo\`{e}ve, M.: Probability Theory II.
Springer-Verlag New York, 4th edition (1978). 

\bibitem{LordPowellShardlow2014}
Lord, G.J., Powell, C.E., Shardlow, T.: An Introduction to Computational Stochastic PDEs.
Cambridge University Press (2014). 

\bibitem{Fishman1996}
Fishman, G.: Monte Carlo: Concepts, Algorithms, and Applications.
Springer-Verlag New York (1996). 

\bibitem{BrennerScott2008}
Brenner, S.C., Scott, L.R.: The Mathematical Theory of Finite Element Methods.
Springer-Verlag New York, 3rd edition (2008).

%

\bibitem{BarraultMadayNguyenPatera2004}
Barrault, M., Maday, Y., Nguyen, N.C., Patera, A.T.: An ‘empirical interpolation’ method: application to efficient reduced-basis discretization of partial differential equations.
Comptes Rendus Mathematique, (2004). \url{doi:10.1016/j.crma.2004.08.006}

\bibitem{HuynhRozzaSenPatera2007}
Huynh, D.B.P., Rozza, G., Sen, S., Patera, A.T.:A successive constraint linear optimization method for lower bounds of parametric coercivity and inf–sup stability constants.
Comptes Rendus Mathematique (2007). \url{doi:10.1016/j.crma.2007.09.019}

\bibitem{Haadonk2017}
Haasdonk, B.: Reduced Basis Methods for Parametrized PDEs-- A Tutorial Introduction for Stationary and Instationary Problems. 
In Benner, P., Ohlberger, M., Cohen, A., Willcox, K.: Reduction and Approximation: Theory and Algorithms, chapter 2.
Society for Industrial and Applied Mathematics (2017). \url{doi:10.1137/1.9781611974829.ch2}



\bibitem{ChenQuarteroniRozza2017}
Chen, P., Quarteroni, A., Rozza, G.: Reduced Basis Methods for Uncertainty Quantification.
SIAM/ASA Journal on Uncertainty Quantification (2017). \url{doi:10.1137/151004550}

\bibitem{PrudhommeRovasVeroy2001}
Prud'homme, C., Rovas, D.V., Veroy, K., Machiels, L., Maday, Y., Patera, A.T., Turinici, G.: Reliable Real-Time Solution of Parametrized Partial Differential Equations: Reduced-Basis Output Bound Methods. 
ASME. J. Fluids Eng. (2001). \url{doi:10.1115/1.1448332} 

\bibitem{GubischVolkwein2017}
Gubisch, M., Volkwein, S.: Proper Orthogonal Decomposition for Linear-Quadratic Optimal Control. 
In Benner, P., Ohlberger, M., Cohen, A., Willcox, K.: Reduction and Approximation: Theory and Algorithms, chapter 1.
Society for Industrial and Applied Mathematics (2017). \url{doi:10.1137/1.9781611974829.ch1}


\end{thebibliography}
\end{document}